\documentclass[11pt]{article}

\usepackage{fullpage}
\usepackage{setspace}
\usepackage{amsmath}
\usepackage{amssymb}
\usepackage{amsthm}
\usepackage{latexsym}
\usepackage{color}
\usepackage{graphicx}
\usepackage{color}
\usepackage[T1]{fontenc}
\usepackage[english]{babel}
\usepackage{geometry}
\usepackage{enumerate}

\DeclareSymbolFont{calletters}{OMS}{cmsy}{m}{n}
\DeclareSymbolFontAlphabet{\mathcal}{calletters}

\def\be{\begin{equation}}
\def\ee{\end{equation}}

\def\b*{\begin{equation*}}
\def\e*{\end{equation*}}

\newtheorem{Theorem}{Theorem}[part]
\newtheorem{Definition}{Definition}[part]
\newtheorem{Proposition}{Proposition}[part]

\newtheorem{Assumption}{Assumption}[part]
\newtheorem{Lemma}{Lemma}[part]

\newtheorem{Remark}{Remark}[part]
\newtheorem{Example}{Example}[part]

\makeatletter \@addtoreset{equation}{section}

\@addtoreset{Definition}{section}

\@addtoreset{Theorem}{section}

\@addtoreset{Proposition}{section}

\@addtoreset{Property}{section}

\@addtoreset{Assumption}{section}

\@addtoreset{Corollary}{section}

\@addtoreset{Lemma}{section}

\@addtoreset{Remark}{section}

\@addtoreset{Example}{section}

\newcommand{\rmi}{{\rm (i)$\>\>$}}
\newcommand{\rmii}{{\rm (ii)$\>\>$}}
\newcommand{\rmiii}{{\rm (iii)$\>\>$}}

\newcommand{\No}[1]{\left\|#1\right\|}     
\newcommand{\abs}[1]{\left|#1\right|}     

\def \E{\mathbb{E}}
\def \F{\mathbb{F}}
\def \H{\mathbb{H}}
\def \I{\mathbb{I}}
\def \L{\mathbb{L}}
\def \M{\mathbb{M}}

\def \P{\mathbb{P}}
\def \Q{\mathbb{Q}}
\def \R{\mathbb{R}}
\def \S{\mathbb{S}}

\def\Ec{{\cal E}}
\def\Fc{{\cal F}}

\def\Mc{{\cal M}}

\def\Tc{{\cal T}}

\def\Yc{{\cal Y}}
\def\Zc{{\cal Z}}

\def\Yb{{\bar Y}}
\def\Zb{{\bar Z}}

\def \Sum{\displaystyle\sum}
\def \Prod{\displaystyle\prod}

\def\esup{{\rm ess \, sup}}

\def\x{\times}

\def\={\;=\;}
\def\.{\;.}

\def\eps{\varepsilon}

\def\reff#1{{\rm(\ref{#1})}}

\def\1{{\bf 1}}
\def \ep{\hbox{ }\hfill{ ${\cal t}$~\hspace{-5.1mm}~${\cal u}$   } }

\def \proof{{\noindent \bf Proof. }}
\def \ep{\hbox{ }\hfill$\Box$}

\def\Om{\Omega}
\def\om{\omega}

\def\normeL2#1{\left\|{#1}\right\|_{L^2}}

\def\NoH#1#2#3{\No{#1}^{#2}_{{\mathbb{H}^{#2,#3}}}}
\def\NoHu#1#2#3{\No{#1}^{#2}_{{\mathbb{H}^{#2,#3}_{1}}}}
\def\NoM#1#2#3{\No{#1}^{#2}_{{\mathbb{M}^{#2,#3}}}}
\def\NoI#1#2#3{\No{#1}^{#2}_{{\mathbb{I}^{#2,#3}}}}
\def\NoS#1#2{\No{#1}^{#2}_{{\mathbb{S}^{#2}}}}
\def\NoL#1#2{\No{#1}^{#2}_{{\mathbb{L}^{#2}}}}

\title{A unified approach to {\sl a priori} estimates for supersolutions of BSDEs in general filtrations}

\author{Bruno Bouchard\footnote{Universit\'e Paris-Dauphine, PSL Research University, CNRS, UMR [7534], CEREMADE, 75016 PARIS, FRANCE. ANR Liquirisk. bouchard@ceremade.dauphine.fr} 
\and Dylan Possama\"i\footnote{Universit\'e Paris-Dauphine, PSL Research University, CNRS, UMR [7534], CEREMADE, 75016 PARIS, FRANCE. possamai@ceremade.dauphine.fr} 
\and Xiaolu Tan\footnote{Universit\'e Paris-Dauphine, PSL Research University, CNRS, UMR [7534], CEREMADE, 75016 PARIS, FRANCE. tan@ceremade.dauphine.fr, the author gratefully acknowledges the financial support of the ERC 321111 Rofirm, the ANR Isotace, and the Chairs Financial Risks (Risk Foundation, sponsored by Soci\'et\'e G\'en\'erale) and Finance and Sustainable Development (IEF sponsored by EDF and CA).}
\and
	Chao {\sc Zhou}\footnote{Department of Mathematics, National University of Singapore, Singapore, matzc@nus.edu.sg. Research supported by NUS Grants R-146-000-179-133 and R-146-000-219-112.}}
             \date{\today}

\begin{document}

 \maketitle

 \begin{abstract}

	We provide a unified approach to {\sl a priori} estimates for supersolutions of BSDEs in general filtrations, which may not be quasi left-continuous. Unlike the previous related approaches in simpler settings, our results do not only rely on a simple application of It\^o's formula and classical estimates, but use crucially appropriate generalizations of deep estimates for supermartingales obtained by Meyer. As an example of application, we prove that reflected BSDEs are well-posed in a general framework which has not been covered so far in the existing literature. 


%
%
\vspace{3mm}

\noindent{\bf Key words:} Backward SDE, supersolution, Doob-Meyer decomposition, reflected backward SDE. 

\vspace{3mm}

\noindent{\bf MSC Classification (2010):} 60H99

\end{abstract}

\section{Introduction}
Supersolutions of backward stochastic differential equations (BSDEs from now on) were introduced by El Karoui et al. in their seminal paper \cite{ekpq}, in order to  study superhedging strategies in mathematical finance. In the simple context of a filtered probability space $(\Omega,\Fc,\F:=(\Fc_t)_{0\leq t\leq T},\P)$ where $\F$ is the (augmented) natural filtration of a $d$-dimensional Brownian motion $W$, a supersolution of a BSDE with terminal condition $\xi$ and generator $g$ consists in a triple of $\F$-adapted processes $(Y,Z,K)$, living in appropriate spaces, with $K$  predictable non-decreasing, such that
\begin{equation}\label{eq:eqbsde}
Y_t=\xi-\int_t^Tg_s(Y_s,Z_s)ds-\int_t^TZ_s\cdot dW_s+\int_t^TdK_s,\ t\in[0,T],\ \P-{\rm a.s.}
\end{equation}

These objects appeared later to be at the very heart of the study of reflected BSDEs, as introduced in El Karoui et al.~\cite{ekkppq}, and more generally of BSDEs  satisfying  some constraint, see Cvitani\'c, Karatzas and Soner~\cite{CKS} for constraints on  the $Z$-component and Peng et al. \cite{lp,peng,px,px2,px3} for general restrictions.
More recently, supersolutions of BSDEs have been proved to provide the semimartingale decomposition of the so-called second order BSDEs, introduced by Soner, Touzi and Zhang \cite{stz}  and generalized by Possama\"{i}, Tan and Zhou \cite{dxc}, and of  the weak BSDEs studied by Bouchard, \'Elie and R\'eveillac in \cite{ber}.

\vspace{0.5em}
When the generator $g$ is equal to $0$, the process $Y$ defined above is nothing else but a supermartingale, and \reff{eq:eqbsde} is simply its Doob-Meyer decomposition. This was generalized by Peng \cite{peng} using the notion of non-linear supermatingales, see also \cite{BPT15,grigo} and the references therein. 

\vspace{0.5em}
As seen through  the above examples, supersolutions of BSDEs appear quite frequently in the literature, as natural semimartingale decompositions for various stochastic processes, and are often used to study their fine properties.  
Having at hand {\it a priori} estimates on the moments and on the stability of supersolutions is crucial in these contexts.  Unfortunately,  in almost all the previously cited works, with the exception of \cite{dxc}, such estimates have been written in, roughly speaking, the context of a Brownian filtration. This is rather limiting from the point of view of both the theory and the applications, and it has  created a tendency in the recent literature to   reproduce very similar proofs   every time that the context was generalized. 

\vspace{0.5em}
In this paper, we propose a general approach which allows one to consider a quite sufficiently general setting. In particular, we do not assume that the underlying filtration is generated by a Brownian motion. In this case, one needs  to introduce another component in the definition of a supersolution of a BSDE, namely a martingale $M$ that is orthogonal to $W$: 
\begin{equation}\label{eq:eqbsdegen}
Y_t=\xi-\int_t^Tg_s(Y_s,Z_s)ds-\int_t^TZ_s\cdot dW_s-\int_t^TdM_s+\int_t^TdK_s,\ t\in[0,T],\ \P-{\rm a.s.}
\end{equation}
When $K\equiv 0$, such objects were first introduced by El Karoui and Huang \cite{ekh}, and studied more recently by Kruse and Popier \cite{kp} to handle {more }general filtrations, in the context of $\mathbb L^p$-solutions, as in the seminal papers \cite{ekpq,bdm}. Supersolutions in general filtrations play a crucial role for the class of reflected BSDEs studied by Klimsiak \cite{klim3}, which is, as far as we know, the most general reference to date. However, all \cite{ekh}, \cite{klim3}\footnote{{Notice that in \cite{klim3}, when the generator does not depend on $Z$ and $p=1$, there is no need for the quasi left-continuity assumption. But the general case requires it.}}, {\cite{kp}} still impose that the filtration is quasi left-continuous, a property which, for instance, is not satisfied for the second order BSDEs studied in \cite{dxc}. We remind the reader that the filtration $\F$, assumed to satisfy the usual hypotheses, is said to be quasi-left continuous if, for any $\F$-predictable stopping time $\tau$, one has
$$\mathcal F_\tau=\mathcal \Fc_{\tau-}, \text{ where $\mathcal F_{\tau-}:=\sigma\left(A\cap \{t<\tau\},\ A\in\mathcal F_t\right)\vee\mathcal N$},$$
where $\mathcal N$ is the the set of all null sets in $(\Omega, \Fc,\P)$. Intuitively, this means that martingales with respect to $\F$ cannot have predictable times of jumps, and in particular deterministic times of jumps.

\vspace{0.5em}
 To understand the simplifications induced by the quasi left-continuity assumption, let us give a brief sketch of the strategy of proof usually used to obtain estimates, say  in $\mathbb L^2$ for simplicity: 

\vspace{0.5em}
\rmi Apply It\^o's formula to $e^{\alpha \cdot}Y^2$ to obtain
\begin{align}
&e^{\alpha t}Y_t^2+\alpha\int_t^Te^{\alpha s}Y_s^2ds+\int_t^Te^{\alpha s}\No{Z_s}^2ds+\int_t^Te^{\alpha s}d[M]_s \label{intro est}\\
&=e^{\alpha T}\xi^2-2\int_t^Te^{\alpha s}Y_sg_s(Y_s,Z_s)ds-2\int_t^Te^{\alpha s}Y_sZ_s\cdot dW_s-2\int_t^Te^{\alpha s}Y_{s-}dM_s+2\int_t^Te^{\alpha s}Y_{s-}dK_s.\nonumber
\end{align}
\rmii Take expectations on both sides, use classical inequalities (namely Young and Burkholder-Davis-Gundy) and some continuity assumptions on $g$ (usually Lipschitz continuity) to control $\mathbb L^2$-type norms of $(Y,Z,M)$ by the norm of $K$ times a small constant,  when $\alpha>0$ is large enough.

\vspace{0.5em}
\rmiii Use the definition of a supersolution to control the norm of $K$ by the norms of $(Y,Z,M)$, and conclude.

\vspace{0.5em}
What is actually hidden in this reasoning is that, because  the martingale $M$ cannot jump at predictable times when the filtration is  quasi left-continuous,   the bracket $[M,K]$ is identically equal to $0$. This is no longer true for general filtrations, in which case  we have to deal with the term $[M,K]$. It turns out to be   difficult to control, which makes this traditional approach not amenable to filtrations that are not quasi left-continuous. 

\vspace{0.5em}

In the case where we look for estimates in $\L^{2}$, this is not a problem because this bracket term is indeed a (local) martingale, see \cite[proof of Lemma 6]{Schw}\footnote{We thank M.~Schweizer for pointing this argument out to us.}, and taking expectations in both sides (up to localization) is enough. This does not work anymore for $p\ne 2$ because then the formulation corresponding to \reff{intro est} involves a non-linear transformation of this martingale. 
 
\vspace{0.5em}
Of course, such a problem only appears when one considers supersolutions, which is not the case in \cite{ekh} or \cite{klim3}. The main problem in \cite{ekh} is that since they consider BSDEs driven by a general c\`adl\`ag martingale $N$, the generator $g$ is integrated with respect to a Stieljes measure related to $d\langle N\rangle$, and, if the filtration is not quasi-left continuous, then $\langle N\rangle$ may have jumps in general, which prevents the technics in \cite{ekh} to be applied. In \cite{kp} however, the problem comes from the fact that they consider so-called BSDEs with jumps, which adds another martingale in the definition of the BSDE. When the filtration is not quasi-left continuous, this martingale can jump at predictable times, which makes the analysis more difficult.

\vspace{0.5em}
The main aim of our paper is to give a general proof of {\it a priori} estimates and stability for supersolutions of BSDEs in a possibly non-quasi-left-continuous filtration. The proof relies on the following property:
$$\text{"It is sufficient to control the norm of $Y$ to control the norm of $(Y,Z,M,K)$."}$$
This is the philosophy of the estimates of Meyer \cite[Theorem 1]{Meyer} that apply to general super-martingales {(see also the generalization in \cite[Thm 3.1]{leng})}. In Section \ref{sec:1}, we show how it can be generalized to the non-linear context of BSDEs. Namely, Theorem \ref{th:main1} below  provides the extension of \cite[Theorem 1]{Meyer} to a supersolution, while Theorem \ref{th:main2} is a version that applies   to the difference of two supersolutions. Both are valid for supersolutions that are only l\`adl\`ag. In Section \ref{sec:rbsde}, we use these results to provide a well-posedness result for reflected BSDEs with a c\`adl\`ag obstacle. When there is no quasi left-continuity assumption on the filtration, this result is not available in the existing literature.

\vspace{0.5em}

\noindent {\bf Notations.}
	For any $l\in\mathbb N\backslash\{0\}$, we denote the usual inner product of two vectors $(x,y)\in\R^l\times \R^l$ by $x\cdot y$. The Euclidean norm on $\R^l$ is denoted by $\No{\cdot}$, and simplified to $\abs{\cdot}$ when $l=1$. 
Let $T>0$ be fixed and let $(\Om, \Fc, \P)$ be a complete probability space, equipped with a filtration $\F = (\Fc_t)_{0 \le t \le T}$ satisfying the usual conditions, and carrying a standard $d$-dimensional $\F-$Brownian motion $W$. Importantly, we do not assume that the filtration is quasi left-continuous. 
	Given $p > 1$ and $\alpha > 0$, we introduce the classical spaces:
	
	\begin{itemize}
		
		\item $\mathbb L^p$ is the space of $\R$-valued and $\Fc_T$-measurable random variables $\xi$ 
		such that 
		$$\No{\xi}^p_{\mathbb L^p}:=\mathbb E \left[\abs{\xi}^p\right]<+\infty.$$
	
		\item $\mathbb S^p $ (resp. $\S^p_r$) denotes the space of $\R$-valued, $\F$-adapted processes $Y$, with $\P$-${\rm a.s.}$ l\`adl\`ag (resp. c\`adl\`ag) paths, such that
		$$\No{Y}_{\mathbb S^p }^p:=\mathbb E \left[\underset{0\leq s\leq T}{\sup}\abs{Y_t}^p\right
			]<+\infty.$$
	
		\item $\mathbb M^{p,\alpha} $ is the space of $\R$-valued, $\F$-adapted martingales $M$, 
		with $\P$-${\rm a.s.}$ c\`adl\`ag paths, such that $M$ is orthogonal to $W$ and
		$$\No{M}_{\mathbb M^{p,\alpha} }^p:=\mathbb E \left[\left(\int_0^Te^{\alpha s}d[M]_s\right)^{\frac p2}\right]<+\infty.$$
	
		\item $\mathbb H^{p,\alpha}$ (resp.~$\H^{p, \alpha}_1$) is the space of $\R^d$-valued (resp.~$\R$-valued) and  $\F$-predictable processes $Z$ such that
		$$\No{Z}_{\mathbb H^{p,\alpha} }^p:=\mathbb E \left[\left(\int_0^Te^{\alpha s}\No{Z_s}^2ds\right)^{\frac p2}\right]<+\infty.$$
	
		\item $\mathbb I^{p,\alpha} $ {(resp.~$\mathbb I^{p,\alpha}_{+}$,~$\I^{p,\alpha}_r$,~$\I^{p,\alpha}_{+,r}$)} denotes the space of $\R$-valued, $\F$-predictable processes with bounded variations $K$, with $\P$-${\rm a.s.}$ l\`adl\`ag (resp. non-decreasing  l\`adl\`ag, c\`adl\`ag, non-decreasing  c\`adl\`ag) paths, such that
		$$
		\No{K}_{\mathbb I^{p,\alpha} }^p:=\mathbb E \left[\left(\int_0^Te^{\frac\alpha2 s}d{{\rm TV}(K)}_s\right)^{p}\right]<+\infty
		$$
and {$K_{0}=0$}. {In the above ${\rm TV}(K)$ denotes the total variation of $K$. }
		\item We also define $\mathbb S^p_{loc}$ as the collection of processes $Y$   such that,  
		for an increasing sequence of stopping times $(\tau_n)_{n \ge 1}$ satisfying $\P(\lim_{n \to \infty} \tau_n = \infty ) = 1$, the localized process $Y_{\tau_n \wedge \cdot}$ belongs to  $\S^p$ for each $n\ge 1$.  The spaces $\mathbb M^{p,\alpha}_{loc}$, $\mathbb H^{p,\alpha}_{loc}$, $\mathbb H^{p,\alpha}_{1, loc}$, 
		$\mathbb I^{p,\alpha}_{loc}$, and $\mathbb I^{p,\alpha}_{+,loc}$ are defined similarly.

	 	\item {Finally, for $\alpha = 0$, we simplify the notation $\M^p := \M^{p,0}$,
		$\H^{p} := \H^{p, 0}$, $\H^p_1 := \H^{p,0}_1$, $\I^p := \I^{p, 0}$, $\I^p_{+} := \I^{p, 0}_{+}$,
		$\I^p_r := \I^{p, 0}_r$ and $\I^p_{+,r} := \I^{p, 0}_{+,r}$.}

	\end{itemize}
 Note that the above spaces do not depend on the precise value of $\alpha$ as we work on the compact time interval $[0,T]$,    two values of $\alpha$ actually provide equivalent norms. Still, we keep the parameter $\alpha$ which, as usual, will be very helpful for many of our arguments. 
 
 \vspace{0.5em}
 
Given a l\`adl\`ag optional process $X$, such that its right-limit process $X^{+}$  is a semimartingale,   and a locally bounded predictable process $\phi$, we define  the stochastic integral as in \cite{LenglartIto}:
$$
(\phi\star X)_{t}:=		\int_0^t \phi_s dX_s 
		~:=~
		\int_0^t \phi_s dX^+_s - \phi_t (X_{t+} - X_t)
,\;t\ge 0.
$$ 
Moreover, we define $\int_t^T \phi_s dX_s := \int_0^T \phi_s dX_s - \int_0^t \phi_s dX_s$.

\section{{\it A priori} estimates}\label{sec:1}

\vspace{0.5em}
	Let us consider a  BSDE with terminal condition $\xi$ and generator $g: [0,T] \x \Om \x \R \x \R^d \longrightarrow \R$. 
	For ease of notations, we denote $g^0_t( \om) := g_t( \om, 0, 0)$. Although, we will have to differentiate between possible values of $p> 1$, this parameter is fixed from now on. 
	The following standing assumption is assumed throughout this  section.
	
	\begin{Assumption}\label{assum:main}
		\rmi   $\xi \in \L^p$, 
		{$g^0  \in \H^{p}_1$}
		and the process $(t,\omega)\longmapsto g_t(\omega,y,z)$ is $\F$-progressively measurable
		for all $(y,z)\in\R\times\R^d$.

		\vspace{0.5em}

		\rmii There exist $(L_y,L_z)\in \R_+^2$, independent of any variables, s.t. for all $(t, \om, y_1, z_1, y_2, z_2) \in [0,T] \x \Om \x (\R \x \R^d)^2$
		\begin{equation}\label{eq:Lip_g}
		\big| g_t(\om, y_1, z_1) - g_t( \om, y_2, z_2) \big| 
		\le
		L_y  |y_1 - y_2| + L_z\No{z_1 - z_2}.
	\end{equation}

	\end{Assumption}
	
	We recall here the definition of a  supersolution.

	\begin{Definition}
	We say that $(Y, Z, M, K)$ is a  solution $($resp.~local solution$)$ of    
	\begin{equation} \label{eq:BSDE}
		Y_t ~=~ \xi -\int_t^T g_s(Y_s,Z_s)ds -\int_t^T Z_s\cdot dW_s-\int_t^T dM_s + K_T-K_t,
	\end{equation}
	if the above holds for any $t\in[0,T]$, $\P-{\rm a.s.}$, 
	and $(Y, Z, M, K)\in\mathbb S^p \times\mathbb H^p \times\mathbb M^p \times\mathbb I^p$
	$($resp. $(Y, Z, M, K)\in\mathbb S^p_{loc} \times\mathbb H^p_{loc} \times\mathbb M^p_{loc} \times\mathbb I^p_{loc})$.  If moreover $K\in \mathbb I^p_{+}$ $($resp. $\mathbb I^p_{+,loc})$, we say that 
	$(Y,Z,M,K)$ is a supersolution $($resp.~a local supersolution$)$ of \reff{eq:BSDE}.
\end{Definition}

\subsection{Estimates for the solution} \label{sec:11}
Our main result says that one can control $(Z,M,K)$ by controlling the component $Y$ of a solution $(Y,Z,M,K)$.
 We emphasize that the general setting we consider here creates additional difficulties that have not been tackled so far in the literature, and which mainly stems from the fact that it is possible for the processes $K$ and $M$ to jump at the same time, when the filtration is not quasi left-continuous. Therefore, the traditional approach which consists in applying It\^o's formula to $\abs{Y}^p$ to derive the desired estimates fails, as this makes the cross-variation between $M$ and $K$ appear, a term that  has no particular  sign and cannot be controlled easily. Our message here is that, in order to obtain such estimates in a general setting, one should rely on a deeper result  from the general theory of processes, namely the estimates obtained in  Meyer \cite{Meyer} for general supermartingales, a version of which  we recall in the Appendix below, see Lemma \ref{lemma:meyer}.
\vspace{2mm}

The following is an extension to the non-linear context.
	\begin{Theorem}\label{th:main1}
		Let $(Y,Z, M, K)\in \mathbb S^p \times\mathbb H^p \times\mathbb M^p \times\mathbb I^p_{+}$ be a   solution of   \eqref{eq:BSDE}. 
		 Then, for any $\alpha\ge 0$, there is a constant $C^{{\alpha}}_{\ref{th:main1}}$  such that
		\begin{equation*}
			\No{Z}^p_{\mathbb H^{p,\alpha} }
			+ \No{M}^p_{\mathbb M^{p,\alpha} }
			+ \No{K}^p_{\mathbb I^{p,\alpha} }
			\leq
			C^{{\alpha}}_{\ref{th:main1}} \Big( 
				\No{\xi}^p_{\mathbb L^p }+ \No{Y}^p_{\mathbb S^p } + \No{g^0}^p_{\H^{p,\alpha}_1}
			\Big).
		\end{equation*}
	\end{Theorem}
	
Before proving this result, we shall establish more general intermediate estimates, that will also be used to control the difference of solutions in Theorem \ref{th:main2} below. They use the notation
$$
N:=Z\star W+M-K.
$$ 
We start with an easy remark.
\begin{Remark} \label{rem: equiv norme N M K} $(i)$ First note that for any $\ell>0$ and   $(a_i)_{1\leq i\leq n}\subset(0,+\infty)$, 
\begin{align}\label{eq: inega dylan}
 (1\wedge n^{\ell-1} )\Sum_{i=1}^na_i^\ell\leq\left(\Sum_{i=1}^na_i\right)^\ell\leq (1\vee n^{\ell-1})\Sum_{i=1}^na_i^\ell.
\end{align}
Let us now consider a solution $(Y,Z, M, K)\in \mathbb S^p \times\mathbb H^p \times\mathbb M^p \times\mathbb I^p$ of  \eqref{eq:BSDE}. Since $W$ and $M$ are orthogonal, \reff{eq: inega dylan} implies that 
\begin{equation}\label{eq:  Z + M-K le  N le Z + M-K}
(1\wedge 2^{\frac{p}{2}-1}) \left(\No{Z}^p_{\H^{p,\alpha} }+\No{M-K}^p_{\mathbb M^{p,\alpha} }\right)\le \No{N}^p_{\mathbb M^{p,\alpha} } \le (1\vee 2^{\frac{p}{2}-1}) \left(\No{Z}^p_{\H^{p,\alpha} }+\No{M-K}^p_{\mathbb M^{p,\alpha} }\right).
\end{equation}
Moreover, if  $K\in \mathbb I_{+}^{p,\alpha}$  then       the Kunita-Watanabe inequality  leads to 
\begin{equation*}
d[M+Z\star W] \leq   2\left(d[N] +d[K] \right)\le  2\left(d[N] +2K_{-}dK+d[K] \right)\le  2\left(d[N] + dK^{2} \right),
\end{equation*}
so that by \reff{eq: inega dylan} 
\begin{align*}
\No{M+Z\star W}_{\mathbb M^{p,\alpha} }^p&\leq 2^{\frac p2}\mathbb E\left[\left(\int_0^Te^{\alpha s}\left(d[N]_s+dK^2_s\right)\right)^{\frac p2}\right]\\
&\leq 2^{\frac p2}\left(1\vee 2^{\frac p2-1}\right)\left(\No{N}_{\mathbb M^{p,\alpha} }^p+{e^{\frac{\alpha p T}{2}}}\No{K}_{\mathbb I^{p,\alpha} }^p\right).
\end{align*}
Hence, since $W$ and $M$ are orthogonal, we finally have
\begin{align}\label{eq: borne M par K et N}
(1\wedge 2^{\frac{p}{2}-1}) (\No{M}_{\mathbb M^{p,\alpha} }^p+ \No{Z}_{\mathbb H^{p,\alpha} }^p) \le \No{M+Z\star W}_{\mathbb M^{p,\alpha} }^p\le {(2^{\frac p2}\vee 2^{p -1 })}\left(\No{N}_{\mathbb M^{p,\alpha} }^p+{e^{\frac{\alpha p T}{2}}}\No{K}_{\mathbb I^{p,\alpha} }^p\right),
\end{align}
in which the left-hand side inequality remains true even if $K$ is not non-decreasing. 

\vspace{0.5em}
$(ii)$ In the following, we shall also use the standard Young's inequality
\begin{equation}\label{eq: ineq carres}
{ a b\le \beta a^{p} + \frac{ b^{q}}{q(\beta p)^{q/p}}, \mbox{ for $a,b\in \R_+$, $\beta>0$, $p, q > 1$ and  $\frac{1}{p} + \frac{1}{q} = 1$.}}
\end{equation}

{$(iii)$ We also emphasize that for $(Y,Z,M)\in \mathbb S^p \times\mathbb H^p \times\mathbb M^p$, the process
$$\int_0^\cdot e^{p\frac\alpha2 s}\phi_p(Y_{s-})d\left(M+Z\star W\right),$$
is a uniformly integrable martingale, where
\begin{equation}\label{eq: def phi}
 \phi_{p}(y)=|y|^{p-1} {\rm sgn}(y)\1_{y\ne 0},~ \mbox{ for $y\in \R$}.
\end{equation}
Indeed,   Burkholder-Davis-Gundy and H\"older's inequalities imply
\begin{align*}
&\E\left[\underset{0\leq t\leq T}{\sup}\abs{\int_0^te^{p\frac\alpha2 s}\phi_p(Y_{s-})d\left(M+Z\star W\right)_{s}}\right]\\
&\leq C\E\left[\sqrt{\int_0^Te^{p\alpha s}\abs{Y_{s-}}^{2p-2}d[M]_s+\int_0^Te^{p\alpha s}\abs{Y_{s}}^{2p-2}\|Z_s\|^2ds}\right]\\
&\leq Ce^{(p-1)\frac\alpha2T}\No{Y}_{\mathbb S^p}^{p-1}\left(\No{M}_{\mathbb M^{p,\alpha}}+\No{Z}_{\mathbb H^{p,\alpha}}\right),
\end{align*}}
for some $C>0$. 
\end{Remark}  

From now on, we use the generic notation $C$, combined with super- and subscripts, to denote  constants in our estimates that only depend on $L_{y}, L_{z}, p$ and $\alpha$.  If they depend on other parameters, this will be made clear.   Although we do not provide their  expressions explicitly, our proofs are written in such a way that the interested reader can easily keep track of them line after line. 
\vspace{0.5em}

In the following, the inequality \reff{eq: esti K p ge 2} is the crucial one, this is the consequence of Meyer \cite{Meyer}.  
\begin{Lemma} \label{lem: intermediate estimate}  Let $(Y,Z, M, K)\in \mathbb S^p \times\mathbb H^p \times\mathbb M^p \times\mathbb I^p$ be a solution of \eqref{eq:BSDE}. 

\vspace{0.5em} 
$(i)$ If $K\in \mathbb I^p_{+}$, then for all   $\alpha\ge 0$ there exists a constant $C^{\alpha}_{\reff{eq: esti K p ge 2}} $ such that 
\begin{align}
\NoI{K}{p}{\alpha}&\le C^{\alpha}_{\reff{eq: esti K p ge 2}}\left(\NoS{e^{\frac{\alpha}{2} \cdot}Y}{p}
+  \NoH{Z}{p}{\alpha}+ \NoHu{g^{0}}{p}{\alpha}\right).\label{eq: esti K p ge 2}
\end{align}
$(ii)$ If $p\ge 2$, then for all $\eps>0$ {there exists  $\alpha> 0$ and}  a constant  $ C^{\eps,\alpha}_{\reff{eq: esti N p ge 2}}$ such that 
\begin{align}
\hspace{-5mm}\No{Y}^p_{\H^{p,\alpha}_{1} }+\NoM{N}{p}{\alpha}&\le \eps\NoHu{g^{0}}{p}{\alpha} + C^{\eps,\alpha}_{\reff{eq: esti N p ge 2}} \left(\NoL{\xi}{p}
+   \NoL{(e^{\alpha\cdot}Y_{{-}} \star N)_{{T}}}{\frac{p}{2}}\1_{{p>2}}+ \E\left[(e^{\alpha\cdot}Y_{{-}} \star K)_{T}\right]^{+}\1_{{p=2}} \right). 
\label{eq: esti N p ge 2}
\end{align}
$(iii)$ If $p\in (1,2)$, then for all $\eps>0$  {there exists  $\alpha> 0$ and} a constant    $C^{\eps,\alpha}_{\reff{eq: esti N p < 2}}$ such that 
\begin{align}
 \NoM{N}{p}{\alpha}&\le  \eps \NoHu{g^{0}}{p}{\alpha}
+ C^{\eps,\alpha}_{\reff{eq: esti N p < 2}}\left(\NoL{\xi}{p}
+  \NoS{e^{\frac{\alpha}{2}\cdot}Y}{p} +    \E[(e^{p\frac{\alpha}{2}\cdot}\phi_{p}(Y_{{-}})\star K)_{T}]^{+}\right),\label{eq: esti N p < 2}
\end{align}
where $\phi_p$ is defined in \eqref{eq: def phi}.
\end{Lemma}

\proof  {(i)} Let us first prove \reff{eq: esti K p ge 2}. A simple application of It\^o's formula implies that 
$$
e^{\frac\alpha2 \cdot}Y_\cdot-\int_0^\cdot e^{\frac\alpha2 s}\left(g_s(Y_s,Z_s)+\frac\alpha2 Y_s\right)ds
$$ 
is a  supermartingale. Moreover, the non-decreasing process in its Doob-Meyer decomposition is $\int_0^\cdot e^{\frac\alpha2 s}dK_s$. Therefore,   Lemma \ref{lemma:meyer}, Assumption \ref{assum:main} and Jensen's inequality provide 
\begin{align}\label{eq:Kest}
\nonumber&(C^p_{\ref{lemma:meyer}})^{-p}\No{K}_{\mathbb I^{p,\alpha} }^p\leq \No{e^{\frac\alpha2 \cdot}Y-\int_0^\cdot e^{\frac\alpha2 s}\left(g_s(Y_s,Z_s)+\frac\alpha2 Y_s\right)ds}_{\mathbb S^p }^p\\
\nonumber&\leq {(1\vee2^{p-1})} \left(\No{e^{\frac\alpha 2\cdot}Y}^p_{\mathbb S^p }+\E \left[\left(\int_0^Te^{\frac{\alpha}{2}s}\left(\abs{g_s(Y_s,Z_s)}+\frac{\alpha}{2}\abs{Y_s}\right)ds\right)^p\right]\right)\\
&\leq {(1\vee2^{p-1})} \left(\left(1+{(1\vee3^{p-1})}T^{{p}}\left(L_y+\frac\alpha2\right)^p\right)\No{e^{\frac\alpha2\cdot}Y}^p_{\mathbb S^p }+{(1\vee3^{p-1})}\left(L_z^p\No{Z}^p_{\H^{p,\alpha} }+{\NoHu{g^{0}}{p}{\alpha}}\right)  \right),
\end{align}
in which the constant  $C^p_{\ref{lemma:meyer}}$ is   as in Lemma \ref{lemma:meyer}. 

\vspace{0.5em}
 {(ii)} We now turn to \reff{eq: esti N p ge 2}. As usual, we apply It\^o's formula to $e^{\alpha\cdot}Y^2$, see \cite[p.538]{LenglartIto} for the case of l\`adl\`ag processes,  use Assumption \ref{assum:main} and \reff{eq: ineq carres},   to obtain 
\begin{align}\label{eq:eq}
\nonumber&e^{\alpha t}Y_t^2+\left(\alpha-\frac{1}{{\eps}}-2L_y-\frac{L_z^2}{{\eta}}\right)\int_t^Te^{\alpha s}Y_s^2ds+(1-\eta)\int_t^Te^{\alpha s}\No{Z_s}^2ds+\int_t^Te^{\alpha s}d[M-K]_s\\
&\leq e^{\alpha T}\abs{\xi}^2+\eps\int_t^Te^{\alpha s}\abs{g^0_s}^2ds - 2\int_t^Te^{\alpha s}Y_{s-}dN_s,
\end{align}
for any  $(\eps,\eta)\in(0,+\infty)^2$. 
Combined with \reff{eq: inega dylan}, this implies that  
\begin{align*}
C_{1} \No{Y}^p_{\H^{p,\alpha}_{1} }+ C_{2}\No{Z}^p_{\H^{p,\alpha} }+\No{M-K}^p_{\mathbb M^{p,\alpha} }\leq &\ 3^{\frac p2-1}\left(e^{\frac p2\alpha T}\No{\xi}^p_{\mathbb L^p }+\eps^{\frac p2} \No{g^0}^p_{\H^{p,\alpha}_1} \right)\\
 &+  3^{\frac p2-1} 2^{\frac p2} \NoL{(e^{\alpha\cdot}Y_{-} \star N)_{T}}{\frac{p}{2}}\1_{p>2}\\
 &+3^{\frac p2-1} 2^{\frac p2} \E\left[(e^{\alpha\cdot}Y_{{-}} \star K)_{T}\right]^{{+}} \1_{p=2},
\end{align*}
where $C_{1}:=(\alpha-\frac{1}{{\eps}}-2L_y-\frac{L_z^2}{{\eta}})^{\frac p2}$, $C_{2}:=(1-\eta)^{\frac p2}$ and where we have used Remark \ref{rem: equiv norme N M K}(iii) in the case $p=2$. Fix $\alpha>0$ and $\eta \in (0,1)$ such that $C_{1}, C_{2}>0$.  We then deduce  \reff{eq: esti N p ge 2} from the right-hand side of \reff{eq:  Z + M-K le  N le Z + M-K} {for $\alpha$ large enough.}

\vspace{0.5em}
{(iii)} It remains to prove \reff{eq: esti N p < 2}.  Since $p<2$, we can not use the Burkholder-Davis-Gundy inequality with exponent $p/2$ to a martingale involving $M$, as it is only c\`adl\`ag. We then follow the approach proposed recently in \cite{kp}. We first appeal to Lemma \ref{lemma:kp} below:
\begin{align*}
&e^{p\frac\alpha2t} \abs{Y_t}^p+\frac{p(p-1)}{2}\int_t^Te^{p\frac\alpha2s}\abs{\phi_{p-1}(Y_s)}d[  N]_s^c
+\frac{\alpha p}{2}\int_t^Te^{p\frac\alpha 2 s}\abs{Y_s}^pds + A_{t}^{T}\\
&\leq  e^{p\frac\alpha 2T}\abs{\xi}^p+p\int_t^Te^{p\frac\alpha2s}\abs{Y_s}^{p-1}\abs{g_s(Y_s,Z_s)}ds-p\int_t^Te^{p\frac\alpha 2s}\phi_{p}(Y_{s-})dN_s,
\end{align*}
in which 
\begin{align*}
A_{t}^{T}:=\frac{p(p-1)}{2}\Sum_{t<s\leq T}e^{p\frac\alpha 2s}\abs{\Delta  N_s}^2\left(\abs{Y_{s-}}^2\vee\abs{Y_{s-}+\Delta  N_s}^2\right)^{\frac p2-1}{\bf 1}_{\abs{Y_{s-}}\vee\abs{Y_{s-}+\Delta  N_s}\neq 0},
\end{align*}
with $\Delta N:=N_{+}-N_{-}$.  
Recalling  Assumption \ref{assum:main} and using \reff{eq: ineq carres}, this shows  that, for any $\beta,\gamma>0$,
\begin{align*}
&e^{p\frac\alpha2t} \abs{Y_t}^p+\left(\frac{p(p-1)}{2}-\beta\right)\int_t^Te^{p\frac\alpha2s}\abs{\phi_{p-1}(Y_s)}d[  N]_s^c+p\left(\frac{\alpha }{2}-L_y-\frac{pL_z^2}{4\beta}\right)\int_t^Te^{p\frac\alpha 2 s}\abs{Y_s}^pds+A_{t}^{T} \\
&\leq  e^{p\frac\alpha 2T}\abs{\xi}^p+p\int_t^Te^{p\frac\alpha2s}\abs{Y_s}^{p-1}\abs{g^0_s}ds-p\int_t^Te^{p\frac\alpha 2s}\phi_{p}(Y_{s-})dN_s\\
&\leq  e^{p\frac\alpha 2T}\abs{\xi}^p+ {\NoHu{g^{0}}{p}{\alpha}}+\frac{p-1}{p^{\frac{p}{p-1}}}  \underset{0\leq s\leq T}{\sup}\abs{e^{\frac\alpha2s}Y_s}^p-p\int_t^Te^{p\frac\alpha 2s}\phi_{p}(Y_{s-})dN_s.
\end{align*}
Let us take $\alpha\geq 2L_y+pL_z^2/(2\beta)$ with $\beta<p(p-1)/2$. Taking expectations on both sides, we obtain  
\begin{align}
\E\left[\int_t^Te^{p\frac\alpha2s}\abs{\phi_{p-1}(Y_s)}d[  N]_s^c+A_{t}^{T}\right]\le&\ C^{1}_{\reff{eq: esti inter N p < 2}} \NoL{\xi}{p}
+  {C^{2}_{\reff{eq: esti inter N p < 2}}}\NoHu{g^{0}}{p}{\alpha}
+ C^{3}_{\reff{eq: esti inter N p < 2}} \NoS{e^{\frac{\alpha}{2}\cdot}Y}{p}\nonumber\\
&+{C^{4}_{\reff{eq: esti inter N p < 2}}}\;\E[(e^{p\frac{\alpha}{2}\cdot}\phi_{p}(Y_{-})\star K)_{T}]^{{+}}, \label{eq: esti inter N p < 2}
\end{align} 
for some explicit constants $(C^{i}_{\reff{eq: esti inter N p < 2}})_{1\leq i\leq 4}$.
We then argue as in \cite[Step 2, Proof of Proposition 3]{kp}\footnote{As pointed out by a referee, there are some inaccuracies in the proof of Proposition 3 in \cite{kp}, especially their inequality (31), which is only valid for predictable integrands. However, this inequality is never used in our proofs, and we only use similar estimates as those of their Equation (33). } { and use \reff{eq: ineq carres} again} to obtain that
\begin{align*}
&\No{ N}^p_{\mathbb M^{p,\alpha} }\leq \left(\frac{2{C^{2}_{\reff{eq: esti inter N p < 2}}}}{ \eps p}\right)^{\frac{1}{p-1}}(2-p)\No{e^{\frac\alpha 2\cdot}Y}_{\mathbb S^p }^p+ \frac{\eps}{{C^{2}_{\reff{eq: esti inter N p < 2}}}}\; \E\left[\int_0^Te^{p\frac\alpha2s}\abs{\phi_{p-1}(Y_s)}d[  N ]^c_s+A_{{0}}^{T}\right],
\end{align*}
for $\eps>0$.
\ep

\vspace{2mm}

We are now in position to complete the proof of Theorem \ref{th:main1}.

\vspace{2mm}

{\bf Proof of Theorem \ref{th:main1}.} {\bf 1.} We first  assume  that $p\ge 2$.  In the course of this proof, we will have to choose $\alpha>0$ large to apply \reff{eq: esti N p ge 2}. However, since the norms in $\{\NoH{\cdot}{p}{\alpha},\alpha>0\}$ (resp.  $\{\NoHu{\cdot}{p}{\alpha},\alpha>0\}$ and $\{\NoM{\cdot}{p}{\alpha},\alpha>0\}$) are equivalent for different values of $\alpha$, this is enough to prove our general result. We first estimate the last term in \reff{eq: esti N p ge 2}: 
\begin{align*}
\NoL{(e^{\alpha\cdot}Y_{-} \star N)_T}{\frac{p}{2}} 
\le 2^{\frac p 2 -1} \left( \NoL{(e^{\alpha\cdot}Y_- \star (M+Z\star W))_T}{\frac{p}{2}}+\NoL{(e^{\alpha\cdot}Y_-\star K)_T}{\frac{p}{2}} \right),
\end{align*} 
in which,    for any $\delta >0$,
\begin{align*}
\NoL{(e^{\alpha\cdot}Y_-\star K)_T}{\frac{p}{2}}\leq \frac{1}{4\delta}\No{e^{\frac \alpha2\cdot}Y}^p_{\mathbb S^p }+\delta\No{K}_{\mathbb I^{p,\alpha} }^p,
\end{align*}
{recall \reff{eq: ineq carres}}, and  
\begin{align*}
\NoL{(e^{\alpha\cdot}Y_- \star (M+Z\star W))_T}{\frac{p}{2}}&\le C^*_p\E \left[\left(\int_0^Te^{2\alpha s}Y_{s-}^2d[M+Z\star W]_s\right)^{\frac p4}\right]
\\
&\leq  \frac{(C^*_p)^2}{4\delta}\No{e^{\frac \alpha2\cdot}Y}^p_{\mathbb S^p }+\delta\No{M+Z\star W}_{\mathbb M^{p,\alpha} }^p,
\end{align*}
with 
\begin{align}\label{eq: catsup}
C^*_p:=\left(\frac p2\vee \frac{p}{p-2}-1\right)^p{\bf 1}_{p>2}+2{\bf 1}_{p=2},
\end{align}
by {Burkholder}'s inequality, see e.g.~\cite[Theorems 8.6 and 8.7]{osekowski}. Combining  the above inequalities with \reff{eq: borne M par K et N} leads to 
\begin{equation}\label{eq: comb inequa avant final}
\NoL{(e^{\alpha\cdot}Y_- \star N)_T}{\frac{p}{2}}\1_{p>2}+\NoL{(e^{\alpha\cdot}Y_- \star K)_T}{\frac{p}{2}}\1_{p=2} \le  \frac{C_{\reff{eq: comb inequa avant final}}^{1}}{\delta}
\No{e^{\frac \alpha2\cdot}Y}^p_{\mathbb S^p }+\delta C^{2}_{\reff{eq: comb inequa avant final}} \left(    \No{N}_{\mathbb M^{p,\alpha} }^p+ \No{K}_{\mathbb I^{p,\alpha} }^p \right),
\end{equation} 
for some explicit constants $C_{\reff{eq: comb inequa avant final}}^{1} $ and $C^{2}_{\reff{eq: comb inequa avant final}}$ that do not depend on $\delta$. 
By inserting the last inequality in \reff{eq: esti N p ge 2}, for e.g.~$\eps=1$ and $\alpha$ chosen appropriately,    we   obtain for $\delta \in (0,1)$:
\begin{align}\label{eq: borne N tout compris}
(1-\delta C^{4}_{\reff{eq: borne N tout compris}} )\NoM{N}{p}{\alpha}\le   C^{1}_{\reff{eq: borne N tout compris}} \NoL{\xi}{p}+C^{2}_{\reff{eq: borne N tout compris}} \NoHu{g^{0}}{p}{\alpha}
+ \frac{C_{\reff{eq: borne N tout compris}}^{3}}{\delta}
\No{e^{\frac \alpha2\cdot}Y}^p_{\mathbb S^p }+\delta C^{4}_{\reff{eq: borne N tout compris}}  \No{K}_{\mathbb I^{p,\alpha} }^p,
\end{align}
in which the constants $(C^{i}_{\reff{eq: borne N tout compris}})_{i\le 4}$ are explicit and do not depend on $\delta \in (0,1)$. 
In view of the left-hand side of \reff{eq:  Z + M-K le  N le Z + M-K} and  \reff{eq: esti K p ge 2},  \reff{eq: borne N tout compris} provides the required bound on $\No{Z}_{\mathbb H^{p,\alpha} }^p$ by choosing $\delta>0$ small enough, so that we can   then use \reff{eq: esti K p ge 2} again to deduce the corresponding bound on $ \No{K}_{\mathbb I^{p,\alpha} }^p$: 
\begin{equation}\label{eq: finale Z K}
			\No{Z}^p_{\mathbb H^{p,\alpha} }
				+ \No{K}^p_{\mathbb I^{p,\alpha} }
			\leq
			C_{\reff{eq: finale Z K}} \Big( 
				\No{\xi}^p_{\mathbb L^p }+ \No{Y}^p_{\mathbb S^p } + \No{g^0}^p_{\H^{p,\alpha}_1}
			\Big),
		\end{equation}
for some constant $C_{\reff{eq: finale Z K}}$. 
Finally, it remains  to appeal to  \reff{eq: borne M par K et N}, \reff{eq: borne N tout compris} and \reff{eq: finale Z K} to obtain the required bound on $\No{M}_{\mathbb M^{p,\alpha} }^p$ and  conclude the proof in the case $p\ge 2$. 

\vspace{0.5em}
{\bf 2.} We now consider the case $p\in (1,2)$. We argue as above except that we now estimate the last term in \reff{eq: esti N p < 2} by using \reff{eq: ineq carres}: 
\begin{align}
\E[(e^{p\frac{\alpha}{2}\cdot }\phi_{p}(Y_-)\star K)_{{T}}] &\le {\frac{p-1}{\left(\delta p^p\right)^{\frac{1}{p-1}}}} \NoS{e^{\frac{\alpha}{2}\cdot}Y}{p} + {\delta}\NoI{K}{p}{\alpha}.\label{eq: esti last term N p < 2}
\end{align}
The latter combined with \reff{eq:  Z + M-K le  N le Z + M-K}, \reff{eq: borne M par K et N}, \reff{eq: esti K p ge 2} and \reff{eq: esti N p < 2}, as in the end of step 1,  provides the required result after choosing $\delta>0$ small enough.  
\ep

\vspace{0.5em}

	When $(Y, Z, M, K)$ is only a local   solution of   \eqref{eq:BSDE},
	all the arguments above hold  true after a localization.
	Then, using Fatou's Lemma, it follows immediately that $(Y, Z, M, K)$ is a   solution.
	We formulate the following result but omit the proof.

\begin{Proposition}
	Let $(Y, Z, M, K)$ be a local   solution of   \eqref{eq:BSDE}.
	Suppose in addition that $Y \in \S^p$.
	Then, $(Y, Z, M, K)$ is a   solution of   \eqref{eq:BSDE}.
\end{Proposition}

\subsection{Difference of solutions and stability}
In this section, we consider two terminal conditions $\xi^1,\xi^2$,  as well as two generators $g^1$ and $g^2$, satisfying Assumption \ref{assum:main}. We then denote by $(Y^i, Z^i, M^i, K^i)$ $\in \mathbb S^p \times\mathbb H^p \times\mathbb M^p \times\mathbb I^p_{+}$ a  solution of     \eqref{eq:BSDE} with terminal condition $\xi^i$ and generator $g^i$, and set  $N^i:=Z^{i}\star W+M^i-K^i$ , $i=1,2$. For notational simplicity, we also define
$$\delta Y:=Y^1-Y^2,\ \delta Z:=Z^1-Z^2,\ \delta M:=M^1-M^2,\ \delta K:=K^1-K^2,\ \delta N:=N^1-N^2,$$
\vspace{-1.9em}
$$\delta g_t(\omega,y,z):=g^1_t(\omega, y,z)-g^2_t(\omega,y,z),\ \text{for all $(t,\omega, y,z)\in[0,T]\times\Omega\times\R\times\R^d$}.$$

By Assumption \ref{assum:main}, we know that there is an $\R$-valued (resp.~$\R^d$-valued), $\F$-progressively measurable  process $\lambda$ (resp.~$\eta$), with $\abs{\lambda}\leq L_y$ (resp. $\No{\eta}\leq L_z$) such that
$$\delta g_{t}:=g^1_t(Y^1_t,Z^1_t)-g^2_t(Y^2_t,Z_t^2)=\delta g_t(Y_t^1,Z_t^1)+\lambda_t\delta Y_t+\eta_t\cdot\delta Z_t.$$ Then, $(\delta Y,\delta Z,\delta M,\delta K)$ satisfies \reff{eq:BSDE} with   driver $\delta g$ and   terminal condition $\delta \xi$. In particular, we can apply to it the results of Remark \ref{rem: equiv norme N M K} and Lemma \ref{lem: intermediate estimate}. 

The main result of this section, Theorem \ref{th:main2}  below, is in the spirit of Theorem \ref{th:main1}: it suffices to control the norm of $\delta Y$ to control the norms of $\delta Z$ and $\delta (M-K)$. Seemingly, it should just be an application of  Theorem \ref{th:main1} to   $(\delta Y,\delta Z,\delta M,\delta K)$ as it satisfies an equation of the form \reff{eq:BSDE}. However, {it is not the case}:
\begin{itemize}
\item[$(i)$] In Theorem \ref{th:main2}, we will only control $\delta (M-K)$ and not $\delta M$ and $\delta K$ separately. Actually, as shown in Example \ref{ex:c} below, there is no hope to control these two processes separately even in the seemingly benign case where $g^1=g^2=0$.

\item[$(ii)$] Actually, Theorem \ref{th:main2}  can not be an immediate consequence of Theorem \ref{th:main1}, because the process $\delta K$ which appears in the dynamics of $\delta Y$ is no longer non-decreasing, and more importantly because the result of Lemma \ref{lemma:meyer} below  does not hold for quasimartingales (instead of supermartingales). However, it is a direct consequence of the intermediate estimates of Lemma \ref{lem: intermediate estimate}, which explains why they have been isolated. 

\item[$(iii)$] If one has a more precise knowledge of the behavior of the non-decreasing processes $K^1$ and $K^2$, then these estimates can actually be improved. We will make this point more clear when we will treat the special case of reflected BSDEs in Section \ref{sec:rbsde}.
\end{itemize}

Let us now state our result.
\begin{Theorem}\label{th:main2} 
		  For any $\alpha\ge 0$, there is a constant $C^{\alpha}_{\ref{th:main2}}$ such that
		\begin{align*}
			\No{\delta Z}^p_{\mathbb H^{p,\alpha} }+\No{\delta (M-K)}^p_{\mathbb M^{p,\alpha} }
			\leq&
			\    C^{\alpha}_{\ref{th:main2}}\left( 
				\No{\delta \xi}^p_{\mathbb L^p }+ \No{\delta Y}^p_{\mathbb S^p } + \No{\delta Y}^{\frac p2\wedge (p-1)}_{\mathbb S^p } +\No{\delta g(Y^1_{\cdot}, Z^1_{\cdot})}_{\H^{p, \alpha}_1}^p  \right).
\end{align*}
The constant  $C^{\alpha}_{\ref{th:main2}}$ depends on $L_{y}, L_{z}, p$ and $\alpha$, as well as  $(\NoS{Y^{i}}{p},\NoL{\xi^{i}}{p},\NoHu{g^{i}(0,0)}{p}{\alpha})_{i=1,2}$. 
\end{Theorem}

\proof  In this proof, we {take $\alpha$ large enough so as to apply the estimates of Lemma \ref{lem: intermediate estimate}}. The general case is deduced by recalling that the different norms are equivalent for different values of $\alpha$, since $[0,T]$ is compact. 

\vspace{0.5em}
{\bf  1.} We first assume that $p\geq 2$. We apply \reff{eq: esti N p ge 2} to $(\delta Y,\delta Z, \delta M,\delta K)$ and obtain
$$\NoM{\delta N}{p}{\alpha}\le \eps\NoHu{\delta g(Y^1_\cdot,Z^1_\cdot)}{p}{\alpha} + C^{\eps,\alpha}_{\reff{eq: esti N p ge 2}} \left(\NoL{\delta \xi}{p}
+   \NoL{(e^{\alpha\cdot}\delta Y_{{-}} \star \delta N)_{{T}}}{\frac{p}{2}} \right). 
$$
Let us estimate the last term in this inequality. We remind the reader that $\delta N=\delta Z\star W+\delta M-\delta K.$ We first use \reff{eq: inega dylan} to obtain
\begin{align*}
 \NoL{(e^{\alpha\cdot}\delta Y_- \star \delta N)_T}{\frac{p}{2}}&\le 2^{\frac{p}{2}-1}\left(
 \NoL{(e^{\alpha\cdot}\delta Y_- \star (\delta M+\delta Z\star W ))_T}{\frac{p}{2}} +\NoL{(e^{\alpha\cdot}\delta Y_- \star \delta K )_T}{\frac{p}{2}}  \right).
 \end{align*}
 We then apply {Burkholder} inequality and obtain
 \begin{align*}
 \NoL{(e^{\alpha\cdot}\delta Y_- \star (\delta M+\delta Z\star W ))_T}{\frac{p}{2}}
&\le C^{*}_{p}\No{e^{\frac{\alpha}{2}\cdot}\delta Y}_{{\mathbb S}^{p}}^{\frac{p}{2}} \No{\delta M+\delta Z\star W}_{\mathbb H^{p,\alpha}}^{\frac{p}{2}}\\
&\leq 2^{\frac{p-2}{4}}C^*_p\No{e^{\frac{\alpha}{2}\cdot}\delta Y}_{{\mathbb S}^{p}}^{\frac{p}{2}}\sum_{i=1}^2\No{ M^i+ Z^i\star W}_{\mathbb H^{p,\alpha}}^{\frac{p}{2}},
 \end{align*}
 where $C^{*}_{p}$ is as in \reff{eq: catsup}, while  
 \begin{align*}
 \NoL{(e^{\alpha\cdot}\delta Y_- \star \delta K)_T }{\frac{p}{2}}
 &\le  \No{e^{\frac{\alpha}{2}\cdot}\delta Y}_{{\mathbb S}^{p}}^{\frac{p}{2}} \No{\delta K}_{\mathbb I^{p,\alpha}}^{\frac{p}{2}}\leq 2^{\frac{p-1}{2}}\No{e^{\frac{\alpha}{2}\cdot}\delta Y}_{{\mathbb S}^{p}}^{\frac{p}{2}} \sum_{i=1}^2\No{ K^i}_{\mathbb I^{p,\alpha}}^{\frac{p}{2}}.
 \end{align*}
 We can then conclude the proof in the case $p\ge 2$ by using \reff{eq:  Z + M-K le  N le Z + M-K} and the bounds of Theorem \ref{th:main1} {applied to $(Z^{i},M^{i},K^{i})_{i=1,2}$.}
 
 \vspace{0.5em}
 {\bf  2.} We now assume that $p\in (1,2)$ and proceed as above but use \reff{eq: esti N p < 2} in place of  \reff{eq: esti N p ge 2}. Namely, since  
 \begin{align*}
\E\left[(e^{p\frac{\alpha}{2}\cdot}\phi_{p}(\delta Y_-)\star \delta K)_T\right] &\le \No{e^{\frac{\alpha}{2}\cdot}\delta Y}_{{\mathbb S}^{p}}^{p-1} \No{\delta K}_{\mathbb I^{p,\alpha}}\leq \No{e^{\frac{\alpha}{2}\cdot}\delta Y}_{{\mathbb S}^{p}}^{p-1}\sum_{i=1}^2\No{ K^i}_{\mathbb I^{p,\alpha}},
\end{align*}
it suffices to use \reff{eq:  Z + M-K le  N le Z + M-K} and the bound of Theorem \ref{th:main1} applied to $K^{1}$ and $K^{2}$.  
\ep

\section{Application to   reflected BSDEs: a general existence result}\label{sec:rbsde}

 The results of the previous section show that  it suffices to control the norm of $Y$ (resp.~$\delta Y$) in order to control the norm  of $(Z, M, K)$ (resp.~$(\delta Z, \delta M, \delta K)$), given a solution $(Y, Z, M, K)$ (resp.~two solutions $(Y^1, Z^1, M^1, K^1)$ and $(Y^2, Z^2, M^2, K^2)$) of \reff{eq:BSDE}.
In most  examples of applications, we know  how to control the norm of  $Y$ (and $\delta Y$).
This is in particular the case in the context of reflected BSDEs (see e.g.~\cite{ekkppq}), BSDEs with constraints (see e.g.~\cite{CKS, px2}),   2nd order BSDE (see e.g.~\cite{stz, dxc}), weak BSDEs (see \cite{ber}).

\vspace{0.5em}

Let us exemplify this in the context of reflected BSDEs. In particular, the following results extend  Klimsiak \cite{klim1,klim2, klim3} to a  filtration that only satisfies the usual conditions, and may not be quasi left-continuous. For sake of simplicity, we restrict to the case of a c\`adl\`ag obstacle, see \cite{grigo} and the references therein for the additional specific arguments that could be used for irregular obstacles. 
Recall that $\mathbb S^p_{ r}$ (resp.~$\mathbb I^p_{r}$, {$\mathbb I^p_{+,r}$}) denote the set of elements of $\mathbb S^p$ (resp.~$\mathbb I^p$, $\mathbb I^p_{+}$) with  c\`adl\`ag path, $\mathbb P$-{\rm a.s.}

\subsection{A priori estimates for reflected BSDEs}
In this section, we assume that Assumption \ref{assum:main} holds for $\xi$ and $g$.

\begin{Definition}\label{def: reflected} Let $S$ be a c\`adl\`ag process  such that $S^+:= S \vee 0\in \mathbb S^p_r $. 
We say that $(Y,Z,M, K)\in\mathbb S^p_{r} \times\mathbb H^p \times\mathbb M^p \times\mathbb I^p_{{+,r}} $ is a solution of the reflected BSDE with lower obstacle $S$  if  
	\begin{equation} \label{eq:RBSDE}
		Y_t = \xi -\int_t^T g_s(Y_s,Z_s)ds -\int_t^T Z_s\cdot dW_s-\int_t^T dM_s + K_T-K_t,
	\end{equation}
	holds for any $t\in[0,T]$  $\P-{\rm a.s.}$, and if 
	$$\begin{cases}
	\displaystyle Y_t\geq S_t,\ t\in[0,T], \ \P-{\rm a.s.},\\
	\displaystyle \int_0^T\left(Y_{s-}-S_{s-}\right)dK_s=0,\ \P-{\rm a.s.} \ (\text{Skorokhod condition})
	\end{cases}$$

\end{Definition}
In order to provide a first estimate on the component $Y$ of a solution, we use the classical linearization procedure. By Assumption \ref{assum:main}, there exists a $\R$-valued (resp.~$\R^d$-valued), $\F$-progressively measurable (resp. $\F$-predictable)  process $\lambda$ (resp.~$\eta$), with $\abs{\lambda}\leq L_y$ (resp.~$\No{\eta}\leq L_z$) such that
$$
g_s(Y_s,Z_s)=g_s^{0}+\lambda_s Y_s+\eta_s\cdot Z_s,\ s\in[0,T] .
$$
Let us define
\begin{equation}\label{eq: def X}
X_s:=e^{-\int_0^s\lambda_sds},\ \text{and}\ \frac{d\Q}{d\P}:=\mathcal E\left(-\int_0^\cdot \eta_s\cdot dW_s\right)_T, W^\Q:=W_\cdot+\int_0^\cdot\eta_sds,
\end{equation}
in which $\cal E$ denotes the Dol\'eans-Dade exponential. 
Then, by Girsanov theorem, $W^\Q$ is a $\Q$-Brownian motion, $M$ is still a $\Q$-martingale orthogonal to $W^\Q$, and we can re-write  the solution of the reflected BSDE \reff{eq:RBSDE} as 
$$\begin{cases}
\displaystyle X_tY_t=X_T\xi-\int_t^TX_sg^{0}_sds-\int_t^TX_sZ_s\cdot dW^\Q_s-\int_t^TX_{s-}dM_s+\int_t^TX_{s{-}}dK_s,\ t\in[0,T],\\
\displaystyle X_tY_t\geq X_tS_t,\ t\in[0,T],\\
\displaystyle \int_0^TX_{s-}\left(Y_{s{-}}-S_{s-}\right)dK_s=0.
\end{cases}$$
One can now use the link between reflected BSDEs and  optimal stopping problems. The proof is classical so that we omit it, see  \cite{ekkppq} for a proof in a Brownian filtration and for a continuous obstacle, or \cite{lep} for a c\`adl\`ag obstacle, and  \cite{eksaintflour} for more results on optimal stopping. We denote by $\mathcal T_{t,T}$ the set of $[t,T]$-valued $\F$-stopping times, while $\E_{t}$ and $\E^{\Q}_{t}$ stands for the $\Fc_{t}$-conditional expectations under $\P$ and $\Q$. 

\begin{Proposition}\label{prop:optstop}
Let $(Y, Z, M, K)$ and $S$ be as in Definition \ref{def: reflected}.  Then, 
$$
X_tY_t=\underset{\tau\in\mathcal T_{t,T}}{\esup}\ \E^\Q_t\left[-\int_t^\tau X_sg^0_sds+X_\tau S_\tau{\bf 1}_{\tau<T}+X_T\xi{\bf 1}_{\tau=T}\right],
$$
and
$$
Y_t=\underset{\tau\in\mathcal T_{t,T}}{\esup}\ \E _t\left[-\int_t^\tau g_s(Y_s,Z_s)ds+ S_\tau{\bf 1}_{\tau<T}+\xi{\bf 1}_{\tau=T}\right],
$$
for all $t\le T$.
\end{Proposition}
Before continuing, let us introduce the solution $(\mathcal Y,\mathcal Z,\mathcal M)\in\mathbb S^p_{r} \times\mathbb H^p \times\mathbb M^p $ of the following BSDE {(well-posedness is a direct consequence of Theorem \ref{thm: existence Lp} below, see also \cite{kp})}, 
\begin{equation}\label{bsde_prop3.1}
\Yc_t = \xi -\int_t^T g_s(\Yc_s,\Zc_s)ds -\int_t^T \Zc_s\cdot dW_s-\int_t^T d\Mc_s,\ t\in[0,T],\ \P-{\rm a.s.}
\end{equation}
A simple application of the comparison result, which can be proved as in \cite[Proposition 4]{kp}, implies  that
\begin{equation}\label{eq: Y ge Yc}
Y_t\geq \Yc_t,\ t\in[0,T],\ \P-{\rm a.s.}
\end{equation}

Let us first show that Proposition \ref{prop:optstop} is actually enough to control the $Y$-term of a solution (or the $\delta Y$-term of the difference of two solutions) and therefore that Theorems \ref{th:main1} and \ref{th:main2} apply to reflected BSDEs.
\begin{Proposition}\label{prop:ref}
$(i)$ Let $(Y, Z, M, K)$ and  $S$ be as in Definition \ref{def: reflected}. Then,    for any $\alpha\geq 0$,
\begin{align*}
\No{e^{\frac\alpha2\cdot}Y}_{\mathbb S^p }^p\leq & \ C_{\ref{prop:ref}}^{\alpha}\left( \No{\xi}^p_{\L^p }+ \No{e^{L_y\cdot}S^+}^p_{\mathbb S^p }
+ \E \left[\left(\int_0^Te^{L_y s}\abs{g^0_s}ds\right)^{p}\right]\right)+C^{\alpha}_\Yc\No{e^{\frac\alpha2\cdot}\Yc}^p_{\mathbb S^p },
\end{align*}
for some constants $C_{\ref{prop:ref}}^{\alpha}$ and $C^{\alpha}_\Yc$ that only depend on $L_{y}$, $L_{z}$ and $\alpha$. Moreover, if we replace $S^+$ by $S$ in the above, we can take $C^{\alpha}_\Yc=0$.

\vspace{0.5em}
$(ii)$ For $i=1,2$, let $(Y^i, Z^i, M^i, K^i)$ and $S^{i}$ be as in Definition \ref{def: reflected} for a generator $g^i$ satisfying Assumption \ref{assum:main} and terminal condition $\xi^{i}\in \mathbb L^{p}$. Then, for any $\alpha\geq 0$,
\begin{align*}
\No{e^{\frac\alpha2\cdot}\delta Y}_{\mathbb S^p }^p\leq &\ \bar C_{\ref{prop:ref}}^{\alpha}\left( \No{\delta \xi}^p_{\L^p }+ \No{e^{L_y\cdot}\delta S}^p_{\mathbb S^p }+ \E \left[\left(\int_0^Te^{L_y s}\abs{\delta g_s(Y_s^1,Z_s^1)}ds\right)^{p}\right]\right),
\end{align*}
for some constant  $\bar C_{\ref{prop:ref}}^{\alpha}$ that only depends on $L_{y}$, $L_{z}$ and $\alpha$\footnote{We are grateful to Marie-Claire Quenez for indicating us a technical problem in the proof in the first version.}.
\end{Proposition}

\proof
$(i)$ First of all, we recall that
$Y\geq \mathcal Y$ on $[0,T]$. 
Next, we deduce from Proposition \ref{prop:optstop} that,  for any $1<\kappa<p$,
\begin{align*}
&\sup_{t\in[0,T]}\left\{e^{\frac\alpha2t}\abs{Y_t}\right\}\\
&\leq  \sup_{t\in[0,T]}\ e^{\left(L_y+\frac\alpha2\right)t}\E^\Q_t\left[\int_0^T e^{L_ys}\abs{g^0_s}ds+\sup_{s\in[0,T]}\left\{e^{L_ys}S^+_s\right\}+e^{L_yT}\abs{\xi}\right]+\sup_{t\in[0,T]}\left\{e^{\frac\alpha2t}\abs{\Yc_t}\right\}\\
&= \sup_{t\in[0,T]}\ e^{\left(L_y+\frac\alpha2\right)t}\E _t\left[\Ec\left({-}\int_t^T\eta_s\cdot dW_s\right)\left(\int_0^T e^{L_ys}\abs{g^0_s}ds+\sup_{s\in[0,T]}\left\{e^{L_ys}S^+_s\right\}+e^{L_yT}\abs{\xi}\right)\right]\\
&\hspace{0.9em}+\sup_{t\in[0,T]}\left\{e^{\frac\alpha2t}\abs{\Yc_t}\right\}\\
&\leq e^{\left(L_y+\frac\alpha2\right)T+\frac\kappa{2(\kappa-1)}L_z^2T}3^{\frac{\kappa-1}{\kappa}}\sup_{t\in[0,T]} \left(\E _t\left[\left(\int_0^T e^{L_ys}\abs{g^0_s}ds\right)^\kappa+\sup_{s\in[0,T]}\left\{e^{\kappa L_ys}(S^+_s)^\kappa\right\}+e^{\kappa L_yT}\abs{\xi}^\kappa\right]\right)^{\frac 1\kappa}\\
&\hspace{0.9em}+\sup_{t\in[0,T]}\left\{e^{\frac\alpha2t}\abs{\Yc_t}\right\},
\end{align*}
in which it is clear that we could have suppressed the term involving $\Yc$ if we had used $S$ instead of $S^+$.

\vspace{0.5em}
Hence, we deduce from Doob's inequality with the exponent $p/\kappa$ and \eqref{eq: inega dylan} that
\begin{align*}
\No{e^{\frac\alpha2\cdot}Y}_{\mathbb S^p }^p\leq&\ e^{p\left(L_y+\frac\alpha2\right)T+p\frac\kappa{2(\kappa-1)}L_z^2T}6^{p-1}\left(\frac p{p-\kappa}\right)^{\frac p\kappa}\E \left[\left(\int_0^T e^{L_ys}\abs{g^0_s}ds\right)^{p}\right]\\
&+e^{p\left(L_y+\frac\alpha2\right)T+p\frac\kappa{2(\kappa-1)}L_z^2T}6^{p-1}\left(\frac p{p-\kappa}\right)^{\frac p\kappa}\left(\No{e^{L_y\cdot}S^+}_{\mathbb S^p }^p+e^{p L_yT}\No{\xi}^p_{\L^p }\right)+2^{p-1}\No{e^{\frac\alpha2\cdot}\Yc}_{\mathbb S^p }^p,
\end{align*}
where we have to replace the $6^{p-1}$ by $3^{p-1}$ if we use $S$ instead of $S^+$.

\vspace{0.5em}
$(ii)$ We first use a classical argument. We know that there exists an $\R$-valued (resp.~$\R^d$-valued), $\F$-progressively measurable process $\widetilde \lambda$ (resp.~$\widetilde \eta$), with $|\widetilde \lambda|\leq L_y$ (resp.~$\No{\widetilde \eta}\leq L_z$) such that
$$g_s^1(Y_s^1,Z_s^1)-g^2_s(Y_s^2,Z_s^2)=\delta g_s(Y_s^1,Z_s^1)+\widetilde \lambda_s \delta Y_s+\widetilde \eta_s\cdot \delta Z_s,\ s\in[0,T],\ \P-a.s.$$
Therefore, we can define $\tilde \Q\sim \P$ and a bounded positive process $\tilde X$ as in \reff{eq: def X} above such that 
\begin{align*}
\tilde X_{t}{\delta Y_t}&= \E^{\tilde \Q} _t\left[\int_t^{\tau} \tilde X_{s} {\delta g_s(Y^1_s,Z^1_s)}ds+ \tilde X_{\tau}(Y^{1}_{\tau}-Y^{2}_{\tau}) +\int_{t}^{\tau}\tilde X_{s}d(K^{1}-K^{2})_{s}\right],
\end{align*}
for all stopping time $\tau\ge t$.
Set $\tau_{\eps}:=\inf\{s\ge t: Y^{1}_{s}\le S^{1}_{s}+\eps\}\wedge T$.  Clearly, $Y^{1}_{-}\ge  S^{1}_{-}+\eps$ on $[\![ t,\tau_{\eps}]\!]$.  Hence, $K^1_{\tau_{\eps}}-K^1_{t}=0$ by the Skorokhod condition. Moreover, $Y^{1}_{\tau_{\eps}}\le S^{1}_{\tau_{\eps}}+\eps$ on $\{\tau_{\eps}<T\}$. Then, the above leads to
\begin{align*}
\tilde X_{t}{\delta Y_t}&\le \E^{\tilde \Q} _t\left[\int_t^{\tau_{\eps}} \tilde X_{s}\abs{\delta g_s(Y^1_s,Z^1_s)}ds+ \tilde X_{\tau_{\eps}}(Y^{1}_{\tau_{\eps}}-Y^{2}_{\tau_{\eps}}) +0-\int_{t}^{\tau_{\eps}}\tilde X_{s}dK^{2}_{s}\right]
\\
&\leq \underset{\tau\in \Tc_{t,T}}{\esup}\ \E^{\tilde \Q} _t\left[\int_t^\tau \tilde X_{s} \abs{\delta g_s(Y^1_s,Z^1_s)}ds+\tilde X_{\tau} \abs{\delta S_\tau}{\bf 1}_{\tau<T}+\tilde X_{T}\abs{\delta\xi}{\bf 1}_{\tau=T} + \eps \tilde X_{\tau} \right].
\end{align*}
Since the same applies to $Y^{2}-Y^{1}$ in place of $\delta Y=Y^{1}-Y^{2}$, it follows that 
\begin{align*}
\tilde X_{t}\abs{\delta Y_t}&\le \underset{\tau\in \Tc_{t,T}}{\esup}\ \E^{\tilde \Q} _t\left[\int_t^\tau \tilde X_{s} \abs{\delta g_s(Y^1_s,Z^1_s)}ds+ \tilde X_{\tau} \abs{\delta S_\tau}{\bf 1}_{\tau<T}+\tilde X_{T}\abs{\delta\xi}{\bf 1}_{\tau=T} + \eps \tilde X_{\tau} \right].
\end{align*}
It remains to let $\eps$ go to $0$ and then argue exactly as in (i).
\ep

\vspace{0.5em}
We now show that one can actually take advantage of the Skorokhod condition satisfied by the solution of a reflected BSDE to improve the general stability result of Theorem \ref{th:main2}. This result is crucial in order to prove existence of a solution   when $p$ is arbitrary, see the proof of Theorem \ref{rbsde} below. We only provide the result for $p=2$. It could be extended to  $p\in (1,2)$, but this is not important as we can always reduce to $p=2$ by localization, again see the proof of Theorem \ref{rbsde} below.    

\begin{Proposition} \label{eq:rbsde_p2}
For $i=1,2$, let $(Y^i, Z^i, M^i, K^i)$ be as in $(ii)$ of Proposition \ref{prop:ref}. Then, for any $\eps>0$, {there exists   $\alpha$ large enough such that }\begin{align*}
\No{\delta Y}^2_{\mathbb H^{2,\alpha}_{1}}+\No{\delta Z}^2_{\mathbb H^{2,\alpha}}+\No{\delta (M-K)}^2_{\mathbb M^{2,\alpha} }\leq&\   \eps \No{\delta g(Y^1_{\cdot}, Z^1_{\cdot})}_{\H^{2,\alpha}_1}^2+   C_{\ref{eq:rbsde_p2}}^{\alpha}\left(\No{\delta \xi}_{L^2 }^2+ \No{e^{\frac\alpha2\cdot}\delta S}_{\mathbb S^2 } \right),
\end{align*}
for some constant  $  C_{\ref{eq:rbsde_p2}}^{\alpha}$  that only depends on $L_{y}$, $L_{g}$, $\alpha$, $\eps$ and   $(\NoS{Y^{i}}{2},\NoL{\xi^{i}}{2},\NoHu{g^{i}(0,0)}{2}{\alpha})_{i=1,2}$. 
\end{Proposition}

\proof It suffices to apply \reff{eq:  Z + M-K le  N le Z + M-K},  \reff{eq: esti N p ge 2} to $(\delta Y,\delta Z,\delta M,\delta K)$ and use the  Skorokhod condition to deduce the control 
\begin{align*}
\E\left[\int_t^Te^{\alpha s}\delta Y_{s-}d(\delta K_s)\right]&\leq \E\left[\int_t^Te^{\alpha s}\delta S_{s-}d(\delta K_s)\right]
\leq \No{e^{\frac\alpha 2\cdot}\delta S}_{\mathbb S^2 }\No{\delta K}_{\mathbb I^{2,\alpha} },
\end{align*}
in which $\No{\delta K}_{\mathbb I^{2,\alpha}}$ is bounded by Theorem \ref{th:main1} {and Proposition \ref{prop:ref}}. 
\ep

\subsection{Wellposedness of reflected BSDE under general filtration}

	We can finally prove the  existence of a unique solution to the reflected BSDEs \reff{eq:RBSDE}. 
	Our proof is   extremely close to the original one given in \cite{ekkppq}, 
	but relies on the more general estimates given in this paper.

\begin{Theorem}\label{rbsde}
	Let Assumption \ref{assum:main} hold true.
	Then, there is a unique solution $(Y,Z,M,K)$ to
	the reflected BSDE \reff{eq:RBSDE}.  
\end{Theorem}

\proof The uniqueness is an immediate consequence of   Proposition \ref{eq:rbsde_p2}, recall that the Doob-Meyer decomposition of the supermartingale $M-K$ is unique.  We therefore concentrate on the problem of existence. {{In the following, we let $(\Yc,\Zc,\Mc)\in\mathbb S^p_{r} \times\mathbb H^p \times\mathbb M^p$ denote the unique solution to the BSDE (\ref{bsde_prop3.1}) with generator $g$ and terminal condition $\xi$.}} 

\vspace{0.5em}
{\bf 1.} {First, let us consider the case $p = 2$ and prove the existence in 3 steps.}

\vspace{0.5em}
\rmi We assume that the function $g$ does not depend on $(y, z)$, that is, $g_t( \omega, y, z) = g_t( \omega)$. 
	{Using \cite[Theorem 2.12]{klim3},
	there are adapted c\`adl\`ag  processes $(Y, M, K)$ and $Z \in \H^1_{loc}$, such that $M$ is a martingale, $K$ is non-decreasing, and \eqref{eq:RBSDE} holds true.} Moreover, \cite[Corollary 2.8]{klim3} implies 
$$Y_t=\underset{\tau\in\Tc_t }{\esup} ~\E \left[\left.\int_t^\tau g_sds+S_\tau{\bf 1}_{\{\tau<T\}}+\xi{\bf 1}_{\{\tau=T\}}\right|\Fc_t \right],\ t\in[0,T],\ \P-{\rm a.s}.$$
Hence,  
$$ \abs{Y_t}  \leq \E \left[ \left. \abs{\xi}+\int_0^T\abs{g_s}ds+\underset{0 \leq s\leq T}{\sup}S_s^++\underset{0 \leq s\leq T}{\sup}\abs{\Yc_s} \right| \Fc_t \right],$$
{by the above equation combined with \reff{eq: Y ge Yc},} and the Burkholder-Davis-Gundy inequality can be used in conjunction with Theorem \ref{thm: existence Lp} below to deduce that 
$$
 \NoS{Y}{2} \leq C \left(\NoL{\xi}{2}+\No{g}_{\H^2_1} +\NoS{S^+}{2} \right){<\infty},
$$
for some $C>0$. 
{{By Theorem \ref{th:main1}, it follows that $(Z, M, K)\in\mathbb H^2 \times\mathbb M^2 \times\mathbb I^2_{{+,r}} $, and hence $(Y, Z, M, K)$ is the unique solution to \eqref{eq:RBSDE} by Proposition \ref{prop:ref}. 
}}

\vspace{0.5em}

\rmii We now consider  the general case. Given $\alpha>0$, let us consider the space $\mathcal S$ of processes $(Y,Z)$ such that $(Y,Z)\in \mathbb S^2_{r} \x  \H^{2, \alpha} $ and denote
$$
	\No{(Y,Z)}_{2,\alpha}
	:=\No{e^{\alpha \cdot} Y}_{\S^2}+\No{Z}_{\mathbb H^{2,\alpha} }.
$$
Clearly, $(\mathcal S,\No{\cdot}_{2,\alpha})$ is a Banach space. 
Then, existence can be proved by using the classical fixed point argument.
	{Let us set} $(\Yb^0, \Zb^0, \bar M^0,\bar K^0) := (0, 0, 0,0)$ {and define} $(\Yb^n, \Zb^n, \bar M^n,\bar K^n)_{n \ge 1}$ {recursively} as the solution to the following reflected BSDE: 
	\begin{align} \label{eq:rbsde_iter}
		\Yb^{n+1}_t 
		= \xi 
		- \int_{{t}}^{T}  g_{{s}}(\Yb^n_s, \Zb^n_s)  ds
		- \int_{{t}}^{T} \Zb^{n+1}_s\cdot dW_s
		- \int_{{t}}^{T} d \bar M^{n+1}_{{s}}+\int_t^Td\bar K^{n+1},
	\end{align}
	with
	\begin{equation}\label{eq:sko}
	{  \bar Y^{n+1}_t \ge S_t,~ t \in [0,T]}
	~~\mbox{and}~
	\int_0^T\left(\bar Y^{n+1}_{t-}-S_{t-}\right)dK^{n+1}_t=0,\ \P-{\rm a.s.}
\end{equation}
Well-posedness is ensured by the first step above. Indeed, the only thing we have to check is that $g(\bar Y^n_\cdot,\bar Z^n_\cdot)\in\mathbb H_1^2$. However, this is a direct consequence of the fact that $g^0\in\mathbb H_1^2$, that $g$ is uniformly Lipschitz continuous in $(y,z)$, and that $(\bar Y^n,\bar Z^n)\in\mathbb S^2\times\mathbb H^2$ {by  induction}. 

\vspace{0.5em}
\noindent Let us also denote $\bar L^{n}:=\bar M^n-\bar K^n$. Using the estimates of Proposition \ref{eq:rbsde_p2},
	it follows that,  for $\alpha > 0$ large enough, $(\Yb^n, \Zb^n)_{n \ge 1}$ is a Cauchy sequence in $\H^{2, \alpha}_{{1}} \x \H^{2, \alpha}$,
	and hence a Cauchy sequence in $\S^2_{r} \x \H^{2, \alpha}$ by the estimates in Part \rmii of Proposition \ref{prop:ref}.
	Moreover, by \eqref{eq:rbsde_iter}, we have	
	\begin{align}  \label{eq:L_byYZ}
		\bar L^n_t - \bar L^m_t
		=&\ 
		 ( \Yb^n_t - \Yb^m_t) -  (\Yb^n_0 - \Yb^m_0)
		- \int_0^t  (\Zb^n_s - \Zb^m_s) \cdot dW_s
		\nonumber \\
		& -\int_0^t \big( g_{{s}}(\Yb^{n-1}_{{s}}, \Zb^{n-1}_{{s}}) - g_{{s}}( \Yb^{m-1}_{{s}}, \Zb^{m-1}_{{s}})  \big) ds.
	\end{align} 
	It follows that $(\Yb^n, \Zb^n, \bar L^n)_{n \ge 1}$ is a Cauchy sequence in $(\mathcal S\times\mathbb S^2,\No{\cdot}_{2,\alpha}+\No{\cdot}_{\mathbb S^2})$, from which we can pass to the limit and obtain	
		$$\Yb_t 
		= \xi 
		- \int_{{t}}^{T}  g_{{s}}(\Yb_s, \Zb_s)  ds
		- \int_{{t}}^{T} \Zb_s\cdot dW_s
		- \int_{{t}}^{T} d \bar L_{{s}},
$$
as well as $\bar Y\geq S$.

\vspace{0.5em}
\rmiii We now prove that $L$ is a supermartingale having the decomposition $\bar L=:\bar M-\bar K$ where $\bar M$ is orthogonal to $W$ and where the non-decreasing process $K$ satisfies the Skorokhod condition, that is to say that we can pass to the limit in \reff{eq:sko}. First, the fact that $L$ is a c\`adl\`ag supermartingale is immediate from the convergence of $\bar L^n$ in $\mathbb S^2$ and the dominated convergence theorem. The fact that the brackets $[\bar L^n,W]$ converge to $[\bar L,W]$ is clear from, for instance Corollaire $1.9$ in \cite{mem} or the proof of Proposition 2 in \cite{briand}, which proves the orthogonality of $\bar L$ and $W$.

\vspace{0.5em}
Let now $\bar L = \bar M - \bar K$ be its Doob-Meyer decomposition, and let us consider a sequence of stopping times $(\tau_m)_{m \ge 1}$ such that the process 
$\sup_{n \ge 1} (\bar Y^n - S)\1_{[\![0,\tau_{m}[\![}$ is essentially bounded and $\tau_m \longrightarrow \infty$ as $m \longrightarrow \infty$.
Since $\No{\bar Y^n - \bar Y}_{\S^2} + \No{\bar L^n - \bar L}_{\mathbb S^2} \longrightarrow 0$, it follows that
$$
	-\E \left[ \int_0^{\tau_m} ( \bar Y_{t-} -  S_{t-}) d\bar K_t \right] 
	=
	\E \left[ \int_0^{\tau_m} ( \bar Y_{t-} -  S_{t-}) d \bar L_t 
	\right]
	=
	\lim_{n \to \infty} \E \left[ \int_0^{\tau_m} ( \bar Y^n_{t-} - S_{t-}) d \bar L^n_t \right] =0.
$$
Since $\bar K$ is non-decreasing and $\bar Y  \ge S $ on $ [0,T]$, we thus obtain 
$$
\int_0^{\tau_m} ( \bar Y_{t-} - S_{t-}) d \bar K_t = 0, ~\P-{\rm a.s.}
$$
Letting $m \longrightarrow \infty$, we see that the Skorokhod condition \reff{eq:sko} holds true for $\bar K$.

\vspace{0.5em}

{\bf 2.} {Finally, let us consider the general case when $p \in (1, \infty)$.} {{It follows from \reff{eq: Y ge Yc} that one can replace $S$ by $S\vee \Yc$, and therefore reduce to the case where $S\in  \S^p_{r}$ (and not only $S^{+}\in  \S^p_{r}$), which we assume in the following.}}
	In this case, we can define for $n\ge 1$
	$$
		\xi^n := (-n) \vee \xi \wedge n,
		~~~
		S^n := (-n) \vee S \wedge n,
		~~\mbox{and}~~
		{g^n := (-n) \vee g \wedge n,}
	$$
	so that  $(\xi^n, S^n, {g^n( 0,0)})_{n \ge 1} \in \L^2 \x \S^2_{r} \x \H^{2, \alpha}_1$, for any $\alpha\ge 0$. {Thus by Step {\bf 1}, we know that there is a unique solution $(Y^n,Z^n,M^n,K^n)\in\mathbb S^2_{r} \times\mathbb H^2 \times\mathbb M^2 \times\mathbb I^2_{{+,r}}$ to \eqref{eq:RBSDE}.}
	Since   $(\xi^n, S^n, {g^n( 0,0)})_{n \ge 1}$ is a Cauchy sequence in $\L^p \x \S^p_{r} \x \H^{p, \alpha}_1$, 
	the estimates of Proposition \ref{prop:ref} and  Theorem \ref{th:main2}
	show that the sequence $(Y^n, Z^n)_{n \ge 1}$ is also a Cauchy sequence in $\S^p_{r} \x \H^{p, \alpha}$, for any $\alpha\ge 0$.
	Moreover, {by a similar equality as in \eqref{eq:L_byYZ}}, it follows that
	$({L^n := M^n-K^n})_{n \ge 1}$ is also a Cauchy sequence in $\S^p_r$. 
	Using the same arguments as in \rmiii of Step {\bf 1},
	it is easy to check that its limit  is a solution to  \eqref{eq:RBSDE}.
\ep

\section{Side remark}

Note that the existence and uniqueness of 
a solution $(Y, Z, M)\in\mathbb S^p_{r} \times\mathbb H^p \times\mathbb M^p  $ to the BSDE 
	\begin{equation}\label{eq: Bsde Lp}
		Y_t = \xi -\int_t^T g_s(Y_s,Z_s)ds -\int_t^T Z_s\cdot dW_s-\int_t^T dM_s   
	\end{equation}
follows from the same arguments as the one used in Section \ref{sec:rbsde} whenever $\xi \in \mathbb L^{p}$ and $g^{0}\in \mathbb H^{p}_{1}$. Indeed, we can bound the component $Y$ of a solution as in the proof of Proposition \ref{prop:ref} by 
$$
X_tY_t=\E^\Q_t\left[-\int_t^T X_sg^0_sds+X_T\xi \right], ~\P-\mbox{a.s.,}
$$
in which $X$ and $\Q$ are defined as in \reff{eq: def X}. The difference of the $Y$-components  of two solutions can be bounded similarly. Then, it suffices to apply the same fixed point argument as in the proof of Theorem \ref{rbsde}.

\begin{Theorem}\label{thm: existence Lp} Let   Assumption \ref{assum:main} hold true. Then, 
\reff{eq: Bsde Lp} admits a unique solution in $\mathbb S^p_{r} \times\mathbb H^p \times\mathbb M^p$. Moreover, for all $\alpha\ge 0$, there exists a constant $C^{\alpha}_{\ref{thm: existence Lp} }$ that depends only on $L_{y},L_{z}$ and $\alpha$, such that 
\begin{align*}
\No{e^{\frac\alpha2\cdot}Y}_{\mathbb S^p }^p\leq & \ C^{\alpha}_{\ref{thm: existence Lp} }\left( \No{\xi}^p_{\L^p } 
+ \E \left[\left(\int_0^Te^{L_y s}\abs{g^0_s}ds\right)^{p}\right]\right).
\end{align*}
\end{Theorem}

\begin{Remark}\label{rem: diff Lp BSDE} Recalling that the difference of two solutions, with different terminal conditions and generators, is still a solution to a BSDE, the bound of Theorem \ref{thm: existence Lp} applies. Using similar notations as above, we have: 
\begin{align*}
\No{e^{\frac\alpha2\cdot}\delta Y}_{\mathbb S^p }^p\leq & \ \bar C^{\alpha}_{\ref{thm: existence Lp} }\left( \No{\delta \xi}^p_{\L^p } 
+ \E \left[\left(\int_0^Te^{L_y s}\abs{\delta g_s(Y^{1}_{s},Z^{1}_{s})}ds\right)^{p}\right]\right).
\end{align*}
\end{Remark}

\appendix

\section{Appendix}

	{Let us consider a} strong supermartingale $X \in \S^p$ on $[0,T]$. {Then}, its paths are almost surely l\`adl\`ag, and {it} admits the (unique) Doob-Meyer decomposition (see e.g. Mertens \cite{Mertens}) :
	\begin{equation} \label{eq:DoobMeyerMertens}
		X_t = X_0 + M_t - A_t - I_t,
	\end{equation}
	where $M$ is a right-continuous martingale, with $M_0=0$, $A$ is a predictable non-decreasing right-continuous process with $A_0 = 0$,
	and $I$ is a predictable non-decreasing left-continuous process with $I_0 = 0$.
\vspace{2mm}

{The following extends \cite[Thm 3.1]{leng}, which is the key ingredient of our main results, Theorem \ref{th:main1} and Theorem \ref{th:main2}.}
\begin{Lemma}\label{lemma:meyer}
		For every constant $p > 1$, there is some constant $C^{p}_{\ref{lemma:meyer}}  > 0$ such that, 
		for all strong supermartingale $X \in \S^p$ $($with the decomposition \eqref{eq:DoobMeyerMertens}$)$, one has
		$$
			\No{ A}_{\mathbb I^p } 
			+ 
			\No{ I}_{\mathbb I^p } 
			\le
			C^{p}_{\ref{lemma:meyer}} \No{X}_{\mathbb S^p }. 
		$$
	\end{Lemma}
	\proof \rmi We suppose in addition that $X$ is right-continuous so that $I \equiv 0$.
	Denote $X^* := \sup_{0 \le t \le T} |X_t| $ and
	define $\tilde X$ as the right-continuous version of the martingale 
	$\E \big[ X^* \big| \Fc_t \big]$. Then $\widehat X : = X +  \tilde X$ is a non-negative right-continuous supermartingale on $[0,T]$,
	with the Doob-Meyer decomposition
	$$
		\widehat X_t = X_0+( \tilde X_t + M_t ) - A_t.
	$$
	Setting $\widehat X_t  := \widehat X_T$ for $ t \in [T, T+1)$ and $\widehat X_t  = 0$ for $t \in [T+1, \infty)$,
	then $\widehat X$ is in fact a right-continuous potential
	(recall that a potential is a non-negative right-continuous supermartingale $\widehat X$ on $[0, \infty)$ such that $\lim_{t \to \infty} \E[\widehat X_t] = 0$).
	Using {Meyer \cite[Thm 1]{Meyer}} (see also \cite[Thm 3.1]{leng}), there is some constant $C'_p$ such that
	$$
		\No{ A}_{\mathbb I^p } 
		\le
		C'_p \No{\widehat X}_{\mathbb S^p }.
	$$
	By the definition of $\widehat X$ and Doob's martingale inequality, we get
	$$
		\No{ A}_{\mathbb I^p } 
		\le
		C'_p \No{\widehat X}_{\mathbb S^p }
		\le
		C'_p \left(1+\frac{p}{p-1}\right)\No{X}_{\mathbb S^p }.
	$$

\noindent \rmii We now consider the case when $X$ is not necessary right-continuous.
	By Mertens \cite{Mertens}, we know that the process $X+ I$ is a right-continuous strong supermartingale,
	and {that} the left-continuous process $I$ is obtained as the limit of an increasing sequence 
	$I_{\cdot} := \lim_{\eps \to 0} \lim_{n \to \infty} I^{\eps, n}_{\cdot}$, where $I^{\eps, n}$ is defined by
	$$
		I^{\eps, n}_t := \sum_{k = 1}^{n} \big( X_{\sigma^{\eps, n}_k} - X_{\sigma^{\eps, n}_k+} \big) \1_{\sigma^{\eps, n}_k < t},
	$$
	{in which} $(\sigma^{\eps, n}_k)_{1 \le k \le n}$ is the non-decreasing sequence of stopping times which exhausts the first $n$ jumps  of $X$ such that $X_t - X_{t+} \ge \eps$, i.e.
	$$
		\sigma^{\eps,n}_1 ~:=~ \inf\{t \ge 0 ~: X_t - X_{t+} \ge \eps\},
		~~~
		\sigma^{\eps,n}_{k+1} := \inf\{ t > \sigma^{\eps,n}_k ~: X_t - X_{t+} \ge \eps\}, ~~k=1, \cdots, n-1.
	$$
	In Step \rmiii we will show that $\No{I^{\eps, n}}_{\I^p} \le C_p'' \No{X}_{\S^p }$ for some constant $C''_p > 0$ independent of $\eps$ and $n$, then it follows {from the} monotone convergence theorem that
	$\No{I}_{\I^p} \le C_p'' \No{ X}_{\S^p }$, 
	and hence
	$$
		\No{A}_{\I^p} \le C'_p \left(1+\frac{p}{p-1}\right)\No{X + I}_{\S^p}
		\le C'_p \left(1+\frac{p}{p-1}\right) ( 1 + C''_p) \No{X}_{\S^p}.
	$$

\noindent \rmiii It is now enough to prove that $\No{I^{\eps, n}}_{\I^p} \le C_p'' \No{X}_{\S^p }$ for some constant $C_p''$ independent of $(n, \eps)$.
	Notice that the discrete process $\overline X^n := (X_0, X_{\sigma^{\eps, n}_1}, X_{\sigma^{\eps, n}_1+}, \cdots, X_{\sigma^{\eps, n}_n}, X_{\sigma^{\eps, n}_n+})$ is a discrete time supermartingale.
	By interpolation, we can {turn it into} a right-continuous strong supermartingale on $[0,T]$.
	{Then, using} the results in Step \rmi we obtain that
	$$
		\No{I^{\eps, n}}_{\I^p} \le C_p'' \No{X}_{\S^p },
		~~\mbox{with}~ C_p'' := C'_p \left(1+\frac{p}{p-1}\right).
	$$
	\qed

\begin{Remark}
A careful reading of the proof in \cite{Meyer} shows that the constant can be computed explicitly and is given by
$$
	{ C_p':= \min_{2\leq k < p} \left( p \Prod_{j=2}^k\frac{pj}{p-j} \right)^{\frac{k}{p-1}}  {\bf 1}_{p > 2}
	+ \left(\frac{p^2}{p-1} \right)^{\frac{1}{p-1}} {\bf 1}_{p\in(1,2]}.}
$$
\end{Remark}

	Meyer's result show that for a supermartingale $X$ with Doob-Meyer's decomposition $X = X_0+M - A$, 
	we can control $A$ by $X$;
	the following example shows that given two supermartingales, 
	we cannot control the difference of {the $A$ parts} by the difference of {the} supermartingales.
	
	\begin{Example}\label{ex:c}
		{\rm 
		Let $W$ be a one dimensional Brownian motion. Fix  $\eps > 0$ and  let  $V$ be defined by
		$$
			V_t ~:=~ \sum_{k \ge 0} W_{\tau_k} \1_{[\tau_k, \tau_{k +1})}(t),~
			~\mbox{where}~
			\tau_0 := 0,
			~\tau_{k+1} := \inf \{t \ge \tau_k ~: \abs{W_t - W_{\tau_k}} \ge \eps \}.
		$$
		Notice that $V$ is of finite variation with decomposition $V = V^+ - V^-$ 
		such that $V^+$ and $V^-$ are two non-negative non-decreasing and predictable process.
		Let $X^1 := W - V^+$ and $X^2 := - V^-$, then $\sup_t | {X^{1}_t-X^{2}_{t}} |  = \sup_t |W_t - V_t| \le \eps$,
		but {$ V^+ - V^- = V$} cannot be controlled by $\eps$.
	}	\qed
	\end{Example}

	We finally provide a technical lemma used in the paper. {Recall the definition of $\phi_{p}$ in \reff{eq: def phi} and observe that $\abs{\phi_{p-1}(y)}$ $=$ $|y|^{p-2}\1_{y\ne 0}$. }

\begin{Lemma}\label{lemma:kp}
	Let $X$ be a l\`adl\`ag semimartingale.
	Then for all $p\in(1,2)$ and $\alpha>0$, we have $\P-{\rm a.s.}$ for any $t\in[0,T]$
\begin{align*}
e^{p\frac\alpha2t} \abs{X_t}^p
\leq &\ e^{p\frac\alpha 2T}\abs{X_T}^p-p\int_t^Te^{p\frac\alpha2s}\frac{\alpha}{2}|X_s|^{p}ds-p\int_t^Te^{p\frac\alpha 2s}\phi_{p}(X_{s-})  dX_s \\
&-\frac{p(p-1)}{2}\int_t^Te^{p\frac\alpha2s}\abs{\phi_{p-1}(X_{s})}   d[X]_s^c \\
&-
\frac{p(p-1)}{2}\Sum_{t<s\leq T}e^{p\frac\alpha 2s}\abs{X_{s+} - X_{s-}}^2\left(\abs{X_{s^-}}^2\vee\abs{  X_{s+}  }^2\right)^{\frac p2-1}{\bf 1}_{\abs{X_{s-}}\vee\abs{X_{s+}}\neq 0},
\end{align*}
where we denote $X_{T+} := X_T$.

\end{Lemma}

\proof
This is an immediate consequence of a straightforward adaptation of \cite[Lemmas 7, 8 and 9]{kp},
together with the It\^o's formula for l\`adl\`ag processes in \cite[p.538]{LenglartIto}.
\ep

\end{document}